\documentclass{article}\usepackage{amsmath}
\usepackage{amscd}
\usepackage{latexsym}
\usepackage{amssymb}
\usepackage{amsmath}

\renewcommand{\abstractname}

\input xy \xyoption{all}

\title{Deformations associated with rigid algebras\footnote{Both authors wish to thank the  Mathematisches Forschungsinstitut Oberwolfach and the first wishes also to thank Prof. Fujio
Kubo and Hiroshima University for their hospitalities during part of
the preparation of this paper}}
\author{M. GERSTENHABER and A. GIAQUINTO}

\begin{document}
\maketitle
\newtheorem{theorem}{Theorem}
\newtheorem{corollary}{Corollary}
\newtheorem{lemma}{Lemma} 
\newtheorem{Rem}{Remark}
\newcommand{\G}{\Gamma}
\newcommand{\g}{\gamma}
\vspace{-7mm}
\date{}
{\noindent ${}^1$\textit{Department of Mathematics, University of
Pennsylvania, Philadelphia, PA 19104-6395, \linebreak[0]
U.S.A.\linebreak[0]
email:mgersten@math.upenn.edu}\\
\noindent ${}^2$\textit{Department of Mathematics, Loyola
University of Chicago, Chicago, IL 60626-5385 U.S.A.,
email:tonyg@math.luc.edu}}

\abstractname{}
\begin{abstract}\noindent   The deformations of an infinite dimensional algebra may be controlled not just by its own cohomology but by that of an associated diagram of algebras, since an infinite dimensional algebra may be absolutely rigid in the classical deformation theory for single algebras while  depending essentially on some parameters.  Two examples studied here, the function field of a sphere with four marked points and the first Weyl algebra, show, however, that the existence of these parameters may be made evident by the cohomology of a diagram (presheaf) of algebras constructed from the original.  The Cohomology Comparison Theorem asserts, on the other hand, that the cohomology and deformation theory of a diagram of algebras is always the same as that of a single, but generally rather large, algebra constructed from the diagram.
\end{abstract}

\section{Introduction} The deformation theory of infinite dimensional algebras differs from that of finite dimensional ones in several significant ways.  One is that an infinite dimensional algebra $A$ may be absolutely rigid in the classical theory \cite{G:DefI}, i.e., have $H^2(A,A) = 0$, but nevertheless depend essentially on one or more parameters.  This is intrinsic to the theory but raises the question of how to capture, by homological methods, the variability of an infinite dimensional algebra embedded in a family of continuously varying algebras when the algebra itself has no infinitesimal deformations.   Two examples presented here, the function field of the sphere with four marked points and the first Weyl algebra, suggest that while there may be no classical infinitesimal deformations of a given algebra, a suitable diagram (presheaf) of algebras associated with the original will have non-trivial cohomology and will deform as the algebra varies. This passage from a single algebra to a diagram of algebras, whose theory will be reviewed in \S 3, is  analogous to the covering of a manifold by coordinate patches.  Although the constructions in our two examples are relatively
natural, our present understanding of how to construct a 
diagram of algebras associated to a single algebra is incomplete. More elaborate diagrams might give more information about a single algebra (although not, we believe, in the examples here), but while there is surely a point at which nothing more can be extracted we have no criterion for when it is reached.

Another important difference between finite and infinite dimensional algebras is their behavior under specialization. For finite dimensional algebras $A$, the dimensions of the cohomology groups $H^n(A,A)$ are upper semicontinuous functions of any parameters on which $A$ may depend but an infinite dimensional algebra may specialize to a rigid one. At any specialization there is a certain subgroup of ``fragile'' classes which vanish, although simultaneously new classes might appear.

As mentioned, the first algebra of our two examples is that of the sphere with four marked points or ``punctures''. With $\mathbb{C}$ as coefficient ring, if the four marked points on the RIemann sphere are $\infty, 0, 1, \lambda$ then the algebra is $\mathbb{C}[x, 1/x, 1/(x-1), 1/(x-\lambda)]$. The second algebra is the (first) Weyl algebra $W_1 = k\{x,y\}/(xy-yx-1)$, i.e., the non-commutative polynomial ring in two variables $x, y$ with the relation $xy-yx = 1$. Although this has no obvious parameters, one can introduce one to get the $q$-Weyl algebra $W_q = k\{x,y\}/(qxy-yx-1)$.  The algebra of the four-punctured sphere and the Weyl algebra have the common property that denoting either by $A$ one has $H^n(A,A) = 0, n> 0$.  In particular, both are absolutely rigid, notwithstanding that the first depends non-trivially on the parameter $\lambda$ and that one can introduce a parameter into the second.  There is, however, a significant difference between the examples beyond the fact that the first is commutative and the second is not.  The moduli space of the sphere with four marked points is analytically well-behaved  at every $\lambda$, but the moduli space for the Weyl algebra  is not.  This might already be anticipated from the fact $H^*(W_q, W_q)$ is not trivial cf. \cite{GerstGiaq:WeylCohom}.

\section {Some background}  The number of complex parameters or ``moduli'' on which a compact
Riemann surface of genus $g$ depends was known by Riemann: it is zero
for $g=0$ (all Riemann spheres are analytically isomorphic), one for
$g=1$, and $3g-3$ for all $g > 1$. This is also the dimension of the
space of non-singular quadratic differentials, which Teichm\"uller \cite{Teich:Extremale}
identified with the space of infinitesimal deformations, and which he used to produce global extremal quasiconformal deformations. However, Teichm\"uller's concept of an infinitesimal deformation had  no obvious extension to higher dimensions. A major advance due to Froelicher and Nijenhuis \cite{FN}  consisted in showing that if $X$ is a complex analytic manifold of any dimension, and $T$ its sheaf of germs of holomorphic tangent vectors, then the space of infinitesimal
deformations of $X$ can be identified with $H^1(X, T)$. This led to the exhaustive study of analytic deformation theory by K. Kodaira and D. C. Spencer; for references  cf. \cite{Kodaira:CxMfds}.  

A useful
intuitive description of the infinitesimal  deformations of $X$  is due to Spencer: Consider $X$ as `sewn together' from coordinate patches $X_i$; deform $X$ by unstitching them and
letting them slide over each other before resewing. In every overlap
$X_i\cap X_j$, the derivative of the motion of $X_j$ relative to
$X_i$ defines a tangent vector field, and these derivatives give the element of
$H^1(X, T)$ which is the infinitesimal of the deformation. One should view $H^1(X, T)$ as having degree 2, and  more generally $H^q(X, \bigwedge^pT)$
as having degree $p+q$; the dimensions assigned to infinitesimals will then agree with those in the algebraic case. The space of infinitesimal deformations of $X$ as a complex manifold does not exhaust the full space of its infinitesimal deformations, which is $H^2(X, \bigwedge T) = H^0(X,\bigwedge^2T) \oplus  H^1(X,T) \oplus H^2(X, \mathbb{C})$.
The component $H^0(X,\bigwedge^2T)$ is the space of infinitesimal
deformations of $X$ in a non-commutative direction. The
interpretation of the elements of $H^2(X, \mathbb{C})$ as
infinitesimal deformations remains mysterious. (See, however, the remarks at the conclusion of \cite{G:SelfDual}).

While $H^1(X, T)$ is defined by taking a direct limit over
refinements of coverings one can compute it from a single covering
by taking an open covering by Stein manifolds or, in the purely
algebraic case, by affine varieties. To compute $H^*(X, T)$ from a single covering it is essential that the open patches, like  Stein manifolds and affine varieties, have two important properties:  First, the intersection of two such is again of the same kind. Second, from the point of view of classical analytic manifolds or commutative algebraic geometry they have no significant cohomology. Specifically, in the complex case, every holomorphic vector bundle on a Stein manifold is trivial, and in the affine case, cohomology with coefficients in a quasi-coherent sheaf vanishes in positive dimensions. \emph{Note, however, that although the affine patches in  the algebraic case play a role analogous to that of coordinate patches in the analytic case, these affine varieties generally depend on parameters.} That a ``rigid'' algebra may depend on parameters is thus an intrinsic feature of deformation theory.  Once one has a covering by Stein manifolds or affine varieties no further homological information can be extracted by further refinement; it is an analog of this which is one of the aspects missing from the theory presented here.

Spencer's idea can be applied to algebras, keeping in mind
that rings are dual to spaces: in place of the inclusions of
$X_i\cap X_j$ into $X_i$ and $X_j$ we have the morphisms of the
rings of holomorphic functions on $X_i$ and $X_j$ into that of
$X_i\cap X_j$. 

While affine spaces have only trivial cohomology with coefficients in a quasi-coherent sheaf, the cohomology of their function algebras may be intricate. For a regular (smooth) algebra $A$ of characteristic zero,  the Hochschild-Kostant-Rosenberg \cite{HKR} theorem asserts that $H^*(A,A)$ can be identified with the exterior algebra on its derivations; the cohomology with coefficients in an arbitrary module is obtained by tensoring with that module.  Further, in characteristic zero one has the Hodge-type decomposition of $H^*(A,A)$ of Gerstenhaber-Schack, cf \cite{GS:Hodge}. That paper introduced  the homological idempotents in the rational group rings $\mathbb{Q}S_n$ of the symmetric groups (later, but less appropriately, called the Eulerian idempotents). Of these, the first, for every $n$, is `Barr's 
idempotent' \cite{Barr:Harrison} (but so named in \cite{GS:Hodge}).  These homological idempotents partition the Hochschild cochain complex of the algebra into a direct sum of subcomplexes indexed by the positive integers, of which the first, Barr showed, gives Harrison's cohomology for commutative algebras in the case of characteristic zero. In view of  the Hochschild-Kostant-Rosenberg theorem this vanishes in all dimensions greater than one.  All two-dimensional classes can  be identified with skew cocycles and are infinitesimals of deformations to non-commutative algebras.  In characteristic zero, therefore, all regular algebras are absolutely rigid in the commutative theory. Rigidity fails to hold in finite characteristics since certain so-called ``restricted'' $2$-cocycles  then arise as obstructions to the prolongation of derivations to automorphisms, cf. \cite{G:DefIII}, \cite{Skryabin:GroupSchemes}. The class of a restricted two-cocycle is generally not zero and can serve as the infinitesimal of a formally non-trivial  deformation.  However, these deformations  become trivial when the deformation parameter $\hbar$ is replaced by $\hbar^{p^n}$ for some $n$ and so are not ``effective" in the sense that the structure of the algebra has not really changed.  Skryabin \cite{Skryabin:GroupSchemes} proved that when the coefficient field is perfect (a condition which probably can be weakened to being separable generated), that this property conversely characterizes the restricted classes. This suggests that if $A$ is a regular algebra of characteristic $p$ then $H^*(A,A)$  modulo the ideal generated by the restricted classes might again be just the exterior algebra on the derivations of $A$. This would imply that, even in positive characteristics, regular affine algebras have no effective commutative deformations and could serve as coordinate patches if one restricts study to commutative algebras.  

Before discussing the examples, the next section, \S\ref{fragile}, introduces a refinement of the cohomology groups for infinite dimensional algebras, defining the {\em fragile} and the {\em resilient} groups  at a specialization.  
Section \ref{diagram} then
reviews the definition of diagram cohomology and of the
Cohomology Comparison Theorem as given in \cite{GS:Monster}, with simple examples folowing in \S\ref{simple}. However, there have been important
subsequent developments, mainly due to Lowen and Van
den Berg \cite{LowVDBergh:CCT} and Stancu \cite{Stancu:CCT1, Stancu:CCT2} who have shown
that the appropriate setting for the theory is that of derived
categories.

\section{The fragile and resilient cohomology groups}\label{fragile} Unless otherwise mentioned, `the cohomology of an algebra' will mean its cohomology with coefficients in itself. Recall that in general the  cohomology of a deformation of $A$ consists of those classes which can be `lifted' to ones of $A[[\hbar]]$ (with the deformed multiplication) modulo those which lift to coboundaries. When $k$ is a field the cohomology of $A$ has no torsion but that of the deformed algebra may have $\hbar$ torsion classes. These become coboundaries when coefficients are extended to include $\hbar^{-1}$, which is usually tacitly done when $k$ is a field in order that the coefficients continue to form a field. Since the cohomology of the deformed algebra is a subquotient of that of the original,  it follows that the dimensons of the cohomology groups of the deformed algebra can not be larger than those of the original. (If there is a well-defined Euler-Poincar\'e characteristic, however, it is preserved, cf. \cite{GerstGiaq:WeylCohom}.) One might suppose, therefore, that if an algebra depends on a parameter, then the cohomology groups of the algebra obtained at a special value of the parameter can not be smaller than those at a generic value, but that is not always the case for infinite dimensional algebras.  It will be so if the algebra at a generic value can be viewed as a deformation of the specialized one. The Weyl algebra is a specialization of the $q$-Weyl algebra at $q=1$, but the cohomology groups of the latter in positive dimensions are not all zero while those of the Weyl algebra vanish. The Weyl algebra is absolutely rigid and can not be deformed back to the $q$-Weyl algebra in the classical theory.

If  $A$ is an algebra which is finitely generated as a module module over a (commutative, unital) ring $k$ and $\hbar$ a variable, then $A[[\hbar]]$, the ring of all formal series $a(\hbar) =a_0 + a_1\hbar +a_2\hbar^2 + \cdots, \,  a_i \in A$,  is just $A\otimes_k[[\hbar]]$; specializing $\hbar$ to an element of $k$, whenever that is meaningful, therefore simply brings us back to $A$.  However, if $A$ is not finitely generated then $A[[\hbar]]$ is much larger than $A\otimes_k[[\hbar]]$. For if, say, $k$ is a field and if the $a_i$ are linearly independent  then the above $a(\hbar)$ is not in $A\otimes_k[[\hbar]]$.  Now specializing $\hbar$ to an element of $k$, whenever that is defined (which requires, in particular,that $k$ have a suitable topology), may not bring us back to $A$.  There is still a natural isomorphism $H^*(A, A)\otimes k[[\hbar]] \cong H^*(A[[\hbar]], A[[\hbar]])$,  but one has in general little control over the cohomology under specialization. When $A$ is finite dimensional the dimensions of the cohomology groups are upper semicontinuous under specialization but in the infinite dimensional case some classes may vanish. When $k$ is not a field but a ring, additional torsion may appear in the cohomology. Those which become $\hbar$ torsion classes and hence disappear when $\hbar^{-1}$ is adjoined to the coefficients form an ideal in  $H^*(A,A)$; these are the classes which by definition are \emph{fragile} at that specialization.  The quotient of $H^*(A,A)$ by the ideal of fragile classes is the \emph{resilient} cohomology for that specialization.  While the $q$-Weyl algebra does have non-trivial cohomology in dimensions greater than one, those classes are all fragile at the particular specialization $q=1$ since, as mentioned, the Weyl algebra has only trivial cohomology.

\section{Diagram cohomology}\label{diagram} 
A \emph{diagram of algebras} over a small
category $\mathcal{C}$ with objects $i,j,\dots$ is a
contravariant functor $\mathbb{A}$ from $\mathcal{C}$ to the
category of unital associative algebras, i.e., a presheaf of
algebras over $\mathcal{C}$. (One can make the same definition for
diagrams of other kinds of algebras; for the Lie case cf
\cite{GGS:LieDiag}.) For example, the sets in an open covering
$\mathcal{U}$ (closed under taking intersections) of a complex
manifold $\mathcal{M}$ may be viewed as the objects of a category in
which the morphisms are the inclusion maps. That is, if $U, V$ are sets
in $\mathcal{U}$ then there is a unique morphism from $U$ to $V$ if
$U \subset V$ and no morphism otherwise. One
then has a poset (partially ordered set) and any poset may be viewed as a
category in which there is a unique morphism $i \to j$ when $i\le j$
and no morphism otherwise. The functor $\mathbb A$ then assigns to each open set the
algebra of holomorphic functions on that set. In practice one takes
a covering by Stein manifolds and for a smooth algebraic variety one
similarly takes a covering by affine opens since the intersection of
any two is again such and each such set has trivial cohomology; cf
\cite{GS:Monster}. An $\mathbb{A}$-\emph{module} $\mathbb{M}$ is a
contravariant functor from $\mathcal{C}$ to the category of Abelian
groups such that for each $i \in \mathrm{Ob}\,\mathcal{C}$ the group
$\mathbb{M}(i)$ is an $\mathbb{A}(i)$-bimodule and for each morphism
$u:i \to j$ the map $\mathbb{M}(u):\mathbb{M}(j) \to \mathbb{M}(i)$
is a morphism of $\mathbb{A}(j)$-modules where $\mathbb{M}(i)$ is
viewed as an $\mathbb{A}(j)$-module by virtue of the morphism
$\mathbb{A}(u): \mathbb{A}(j) \to \mathbb{A}(i)$.

To define the  cohomology groups $H^n(\mathbb{A}, \mathbb{M})$,
consider first the simplicial complex called the \emph{nerve} of
$\mathcal{C}$. The $0$-simplices of this are just the objects $i$ of
$\mathcal{C}$. For $q > 0$ a non-degenerate $q$-simplex is any
$q$-tuple of composable morphisms $\sigma = (i_0 \to i_1 \to \cdots
\to i_q)$ in which no single morphism $i_r \to i_{r+1}$ is an
identity morphism (although a composite of several of the morphisms
is allowed to be). We will call $i_0$ the \emph{domain} of $\sigma$,
denoted $d\sigma$ and $i_q$ its \emph{codomain}, $c\sigma$. The
$0$-th and $q$-th faces of $\sigma$ are given, respectively, by
$\partial_0\sigma = (i_1\to\cdots \to i_q),\,
\partial_q\sigma = (i_0\to \cdots \to i_{q-1})$. For $0 < r < q$
define $\partial_r\sigma$ by composing the maps $i_{r-1}\to i_r$ and
$i_r \to i_{r+1}$ so  $\partial_r\sigma = (i_0 \to \dots \to i_{r-1}
\to i_{r+1} \to \dots \to i_q)$;  this may be degenerate. Let $C_q(\mathcal{C})$ be the set
of all formal finite linear combinations of non-degenerate
$q$-simplices and set $\partial\sigma = \sum_{r=0}^q
(-1)^r\partial_r\sigma$ omitting any simplices which may be
degenerate. Since the boundary of a degenerate simplex always
vanishes, with this the $C_q$ form a chain complex; it is
isomorphic to the quotient of the complex spanned by all simplices
by the subcomplex spanned by the degenerate simplices. Note that if
$\sigma = (i_0 \to \cdots \to i_q)$ then $\mathbb{M}(d\sigma) =
\mathbb{M}(i_0)$ is a module over $\mathbb{A}(i_q) =
\mathbb{A}(c\sigma)$ by virtue of the composite morphism $i_0 \to
i_q$. 

Now let $C^{p,q}$ be the $k$-module of all functions $\G$ on
$C_q(\mathcal{C})$ which send a $q$-simplex $\sigma$ to a cochain
$\G^{\sigma} \in C^p(\mathbb{A}(c\sigma), \mathbb{M}(d\sigma))$. 
Setting, as before, $\sigma = (i_0 \to \cdots \to i_q)$, those faces
$\partial_r\sigma$ with $0 <r <q$ all have the same domain and
codomain as $\sigma$, but the first and last do not. Write briefly
$\varphi$ for the algebra morphism $\mathbb{A}(i_{q-1} \to i_q):
\mathbb{A}(i_q) \to \mathbb{A}(i_{q-1})$. Then
$\G^{\partial_q\sigma}\varphi$, defined by $(\G^{\partial_q\sigma}\varphi)(a_1,\dots, a_p) =
\G^{\partial_q\sigma}(\varphi a_1,\dots,\varphi a_q), \,a_1,\dots,a_p \in \mathbb{A}(i_q)$ , is again
in $C^p(\mathbb{A}(c\sigma), \mathbb{M}(d\sigma))$. Write $T$ for
$\mathbb{M}(i_0\to i_1):\mathbb{M}(i_1) \to \mathbb{M}(i_0)$. One
then also has $T\G^{\partial_0} \in C^p(\mathbb{A}(c\sigma),
\mathbb{M}(d\sigma))$ and we st
$$\G^{\partial\sigma} = T\G^{\partial_0\sigma}+
\sum_{r=1}^{q-1}(-1)^r\G^{\partial_r\sigma}
+(-1)^q\G^{\partial_q\sigma}\varphi.$$ 
There then exist commuting
coboundaries, the algebraic (Hochschild)
$\delta_{\mathrm{Hoch}}:C^{p,q}\to C^{p+1,q}$ defined by sending $\G^{\sigma}$ to its Hochschild coboundary $\delta\G^{\sigma}$, and the simplicial
$\delta_{\mathrm{simp}}:C^{p,q}\to C^{p,q+1}$ defined by
$(\delta_{\mathrm{simp}}\G)^{\sigma} = \G^{\partial\sigma}$.
Finally, set $C^n(\mathbb{A},\mathbb{M}) = \oplus_{p+q=n}C^{p,q}$
and define the total coboundary $\delta: C^n\to C^{n+1}$ by
$$(\delta\G)^{\sigma} = \delta_{\mathrm{simp}}(\G)^{\sigma}
+ (-1)^{\dim\sigma}\delta_{\mathrm{Hoch}}(\G^{\sigma});$$  the
cohomology groups $H^*(\mathbb{A},\mathbb{M})$ are defined to be those of this complex. 

If the highest dimension of any non-zero cohomology group of any algebra in
the diagram  is $m$ and $d$ is the dimension of the nerve of the underlying category
(the dimension of the largest simplex appearing in it)
 then the diagram as a whole may have non-trivial
cohomology in dimensions up to $m+d$. For
example, denote  by $\mathbf{k}$ the diagram with $\mathbb{A}(i) = k$, the
coefficient ring, and in which every morphism $\mathbb{A}(i\to j)$  is just the
identity.  Then $H^*(\mathbf{k},\mathbf{k})$ is the
simplicial cohomology of the nerve of the underlying category
$\mathcal{C}$ with coefficients in $k$. Here there is no algebra
part; the cohomology is entirely simplicial. (Simplicial cohomology is thus just a special case of 
Hochschild cohomology, cf. \cite{GS:Simplicial}.)

A \emph{deformation of a diagram} of $k$-algebras $\mathbb{A}$ is a
diagram of $k[[\hbar]]$-algebras over the same underlying category
whose reduction modulo $\hbar$ is $\mathbb{A}$, \cite{GS:Monster}.
Because the cohomology of the nerve of
the underlying category $\mathcal{C}$ may not be trivial, unlike the case of a single algebra one can not always identify the
infinitesimal deformations of a diagram $\mathbb{A}$ with
$H^2(\mathbb{A},\mathbb{A})$. In general one must use the cohomology of the {\em asimplicial} subcomplex consisting of those cochains $\G$ where, if the dimension of $\G$ is $n$ and $\sigma$ is an $n$-simplex of the nerve of $\mathcal{C}$ then $\G^{\sigma}$ (which is just an element of $\mathbb{A}(\sigma)$) is rquired to be zero.  However, if all the algebras in the diagram are commutative (as in our first example) or if the geometric realization of the nerve of the underlying category $\mathcal{C}$ is contractible (the case in our second example, but something that could always be accomplished by adjoining a terminator to $\mathcal{C}$) then the problem does not arise, cf \cite{GS:Monster}.  With this refinement, an infinitesimal deformation of $\mathbb{A}$ is again
just the cohomology class of a 2-cocycle $\G$.  The latter
assigns to every $i \in \mathrm{Ob}(\mathcal{C})$ a 2-cocycle
$\G^i$ of $\mathbb{A}(i)$ with coefficients in itself (which we
may interpret as the infinitesimal of a deformation of
$\mathbb{A}(i)$) and assigns to every morphism $i\to j$ a $k$-linear
map $\phi^{ij}:\mathbb{A}(j)\to \mathbb{A}(i)$. Denoting $\phi^{ij}$
for the moment simply by $\phi$, these are connected by the
condition that
$$ \phi(\G^j(a,b))-\G^i(\phi a, \phi b) =
(\delta\G^{ij})(a,b),$$ where $\delta$ is the Hochschild
coboundary. An \emph{integral} of $\G$, if one exists, is a
deformation whose infinitesimal is in the class of $\G$. 

To illustrate the theory we give some simple examples in the following section. The underlying category of a diagram of algebras will always be obvious so we may omit mention of it.

The Cohomology Comparison Theorem (CCT) cf \cite{GS:Monster} asserts that there is a functor which sends a diagram of algebras to a  single associative algebra and a diagrams of modules to a single module which preserves both the cohomology and the deformation theory; the algebra associated to a diagram is its {\em diagram algebra}, and likewise for modules. This allows one to transfer to the cohomology of diagrams of algebras all known properties of the cohomology of a single algebra including, e.g., its Gerstenhaber algebra structure.  The CCT implies, in particular, that the the third cohomology group of a diagram of algebras with coefficients in the diagram itself contains the obstructions to the infinitesimal deformations, so that when this group vanishes every infinitesimal can be integrated.

\section{Some simple diagrams of algebras}\label{simple}  The simplest (non-trivial) diagram of algebras  is just an algebra morphism 
$$\begin{CD}B @>\phi>> A\end{CD}$$  A bimodule over this diagram is a morphism of abelian groups 
$$\begin{CD}N @>T>>M\end{CD}$$
where $N$ is a $B$-bimodule, $N$ an $A$-bimodule, and $T$ a module morphism from $N$ to $M$ with the latter considered as a $B$-bimodule by virtue of the morphism $\phi$.  The original diagram of algebras, which by abuse of notation we may also denote simply by $\phi$ and likewise the module by $T$, is a bimodule over itself.  An $n$-cochain $\G$ of $\phi$ with coefficients in $T$ is a triple  $\G = (\Gamma^{B}, \G^{A}, \G^{\phi})$, the first component of which is a Hochschild $n$-cochain of $B$ with coefficients in $N$, the second an $n$-cochain of $A$ with coefficients in $M$, and the third an ($n\!-\!1$)-cochain of $B$ with coefficients in $M$ considered as a $B$ bimodule. The coboundary is given by  $\delta \G = (\delta \G^{B}, \delta \G^{A}, T\G^{B} - \G^{A}\phi -\delta \G^{\phi})$. This defines the complex $C^*(\phi, T)$;  $\G$ is a cocycle precisely when $\G^{B}$ and $\G^{A}$ are both cocycles and $T\G^{B} - \G^{A}\phi$ is a coboundary. A deformation of $\phi$ is a diagram $\begin{CD}B_{\hbar} @>\phi_{\hbar}>> A_{\hbar}\end{CD}$ where  $B_{\hbar}, A_{\hbar}$ are deformations of $B$ and $A$, respectively and $\phi_{\hbar} = \phi +\hbar\phi_1 +\hbar^2\phi_2 +\dots $ is a morphism between the deformed algebras; here each $k$-linear map $\phi_i:B\to A$ is tacitly extended to be $k[[\hbar]]$ linear. Its infinitesimal is the two-cocycle $\G \in Z^2(\phi, \phi)$ consisting of the first order terms, or more strictly, its cohomology class.  The geometric picture is a morphism $f:X\to Y$ between two topological spaces. If $A$ is the ring of continuous functions on $X$ and $B$ that on $Y$, then $f$ induces a morphism $f^*:B \to A$.   The same idea applies when $X$ and $Y$ are analytic manifolds or varieties defined over some field $k$.

The diagram algebra of a single morphism $\phi:B\to A$ can be written as a direct sum $B+A+A\phi$ where $A\phi$  denotes $A$ considered as a left $A$ and right $B$ module, the latter by the operation of $\phi$. Its elements will be written in the form $a\phi$. With this notation the multiplication is given by $(b, a_1, a_2\phi)(b', a_1', a_2'\phi) =(bb',a_1a_1',a_1a_2'+a_2\phi(b')\phi)$. One can represent the elements $(b, a_1, a_2\phi)$ as  $2 \times 2$ matrices
$$\begin{pmatrix} b & 0 \\ a_2 & a_1
\end{pmatrix}, \quad\text{with}\quad
\begin{pmatrix} b & 0 \\ a_2 & a_1
\end{pmatrix} \begin{pmatrix} b' & 0 \\ a_2' & a_1'
\end{pmatrix} = \begin{pmatrix} bb' & 0 \\ a_2\phi(b')+a_1a_2' & a_1a_1'
\end{pmatrix}.$$
Diagram modules in this case are constructed similarly. The proof of the CCT is not trivial even for this special case.

When deforming a diagram $\phi\!:\!B\to A$, if one  of the algebras, say $B$, is to be held fixed  then one must use the subcomplex of $C^*(\phi, T)$ consisting of those $\G = (0, \G^A, \G^{AB})$ in which the first component is identically zero. (In the geometric situation, suppose that $X$ and the morphism $f$ are deformed but that we require that the space $Y$ to which $X$ maps is held fixed.)  We are then considering only those cocycles $\G^A$ such that $\G^A\phi$ is a coboundary. Since coboundaries are always sent to coboundaries, the cohomology is simply $\ker \phi^*:H^*(A,A)\to H^*(B,A)$.  This is clearly closed under the cup product. It should, we conjecture, also be closed under the graded Lie product because of the relation of that product to obstructions.  Similar considerations hold when $A$ is held fixed.  When both algebras are fixed and we are interested only in how the morphism between them can vary then the infinitesimal is just a derivation of $B$ into $A$; when the morphism is an inclusion and we identify morphisms which differ only by an automorphism of $B$ then the space of infinitesimals consists of the derivations of B into the B bimodule $A/\phi B $, cf \cite{Nij:Subalgs1,NIJ:Subalgs2}.

For the Weyl algebra we will use  a diagram like
$$\begin{CD}\{B @>\phi>>A@<\phi'<<B'\}\end{CD},$$
which will again be denoted by $\mathbb{A}$. A module over this diagram has
the form
$$\begin{CD}\mathbb{M} = \quad  \{N@>T>>M@<T'<<N'\}\end{CD}\qquad$$
where $N$ and $N'$ are, respectively, $B$ and $B'$ bimodules, $M$ is
an $A$ bimodule, $T$ is a $B$ bimodule morphism, $M$ is
considered a $B$ bimodule through the morphism $\phi$, and similarly
for $T'$. An $n$-cochain of $\mathbb{A}$ with coefficients in
$\mathbb{M}$ is a quintuple $\G = (\G^B, \G^{B'},
\G^A, \G^{\phi}, \G^{\phi'})$ as in
$$ \begin{CD}\G =\quad
\{\G^B @>\G^{\phi}>>\G^A
@<\G^{\phi'}<<\G^{B'}\}\end{CD}\qquad.$$ The coboundary is
then 
$$ \begin{CD}\delta\G =\quad
\{\delta\G^B @>T\circ\G^B - \G^A\circ\phi-
\delta\G^{\phi}>>\delta\G^A @<T'\circ\G^{B'} -
\G^A\circ\phi'-
\delta\G^{\phi'}<<\G^{B'}\}\end{CD}\qquad$$ 
where the coboundary
operators inside the braces are the ordinary Hochschild
coboundaries.  The geometric picture is that of a space $X$ with morphisms $f$ and $f'$  (which in general need not be defined on all of $X$) to respective spaces $Y$ and $Y'$. In a very special case, suppose  that we have a diagram 
$$ \quad \begin{CD} B@>f>>A @<g<< C\end{CD} $$ 
in which the second cohomology of every algebra with coefficients in itself vanishes and  where for every derivation $\G^f \in \operatorname{Der}(B,A)$ there is a derivation $\g^A$ of $A$ such that $\G^f(b) = \g^A(fb)$ for all $b\in B$. Then it is easy to check that every $2$-cocycle $\G$ of $\mathbb{A}$ with coefficients itself is cohomologous to one of the form $(0,0,0,0,\G^g)$.

In elementary cases a diagram can not be expected to yield more information than can be directly extracted from the algebras involved. For example, using the theorem of Hochschild-Kostant-Rosenberg it is an easy exercise to show that the cohomology of the diagram
$$k[x] \hookrightarrow k[x,y] \hookleftarrow k[y]$$
can be canonically identified with that of $k[x,y]$ alone.

\section{Deformation of commutative diagrams} There is a simple condition on an infinitesimal deformation of a diagram that insures that it is a deformation `in the direction' of another commutative diagram. It only needs to be stated in the case of a commutative diagram consisting simply of two morphisms and their composition since the generalization is immediate. Consider the diagram

$$\xymatrix{&V\ar[dr]^\beta&\\W\ar[ur]^\alpha\ar[rr]^\theta& &U}$$
 where initially $\theta$  is just $\beta\alpha.$ An infinitesimal deformation assigns, in particular, derivations $\Gamma^\alpha, \Gamma^\beta, \Gamma^\theta$ to the one-simplices $\alpha, \beta, \theta$, respectively. The condition that the infinitesimal respect the commutativity of the diagram is simply that $\beta\Gamma^\alpha + \Gamma^\beta \alpha = \Gamma^\theta$.  If one requires that the composite $\theta$ not change (although the factors $\alpha$ and $\beta$ may) then $\Gamma^\theta = 0$, so we have $\beta\Gamma^\alpha = -\Gamma^\beta \alpha$. This is, in particular, the case when $\beta\alpha=0$ and we require that the composite remain the zero morphism, e.g., when considering the deformation of an exact sequence.

\section{First example: The four-punctured sphere}  Let $k$ be a commutative unital ring and  set $A= k[x, 1/x, 1/(x-1), 1/(x-\lambda)],$ where $\lambda$ is an element of $k$ not equal to 0 or 1. When $k= \mathcal{C}$ this is just the algebra of functions on the Riemann sphere punctured respectively at $\infty, 0, 1, \lambda$.  Since the cohomology of $A$ with coefficients in any $A$ module vanishes in all dimensions greater than 1 it must be that when $A$ is enlarged to $A[[\hbar]]$, where $\hbar$ is a variable, that the dependence on $\lambda$ disappears. This is easy to exhibit  explicitly. Denote $A$ now more explicitly by $A_{\lambda}$. Then inside $A_{\lambda}[[\hbar]]$ we have
$$\frac{x-\lambda}{x-\lambda(1+\hbar)} =
\frac1{1-\hbar\lambda/(x-\lambda)},\quad \mathrm{so}\quad
\frac1{x-\lambda(1+\hbar)} =
\frac1{x-\lambda}\sum_{n=0}^{\infty}\hbar^n\frac{\lambda^n}{(x-\lambda)^n}.$$
Thus $A_{\lambda}[[\hbar]]$ contains $A_{(1+\hbar)\lambda}[[\hbar]]$
but by the same means it is easy to see that the reverse inclusion
also holds, so these rings are identical and hence isomorphic. Note that that the series for $1/(1-\hbar\lambda/(x-\lambda))$ is not contained in $ A_{\lambda}\otimes_kk[[\hbar]]$, but only in the larger ring $A_{\lambda}[[\hbar]]$.  
When $k = \mathbb{C}$ and $x$ is viewed as taking on complex values, the series may be viewed as defining a mapping from the sphere punctured at $\{\infty, 0, 1,  \lambda\}$ to that punctured at $\{\infty, 0, 1, \lambda(1+\hbar)\}$  but there is no value of $\hbar$, however small, for which it can converge everywhere. These surfaces are not analytically isomorphic, but given any neighborhood of $x = \lambda$ one can choose $\hbar$ so small that the series converges everywhere outside that neighborhood. This suggests that analytically the local moduli space is not singular, which is in fact the case. 

The classic deformation theory of a single algebra can not detect the dependence of $A$ as a $k$ algebra on the parameter $\lambda$.  However, if some diagram of $k$ algebras depending only on $\lambda$ and the internal structure of  $A$ deforms non-trivially when $\lambda$ varies then $A$ must vary with $\lambda$ in a non-trivial way inside the family of $k$ algebras, for out of a $k$-isomorphism between algebras one could build an isomorphism of diagrams. 

Suppose that an algebra $A$ depends on a parameter $\lambda$.  A natural attempt at a diagram whose deformation theory captures the dependence would be to present $A$ as a quotient of a free algebra $F$ by a varying ideal $K_{\lambda}$ depending on $\lambda$. 
Unfortunately this presents computational problems. 
Free algebras are generally very large and computing the cohomology of a diagram containing one could be difficult.  We compromise by using a more computable one still showing the dependence of $A$ on $\lambda$ but possibly not capturing all the ways in which $A$ might vary as an algebra over the original ground ring $k$.  

The Riemann sphere punctured at three points does not have any moduli, since any three points $p, q, r$ can be moved by an analytic automorphism to any three others, say $ \infty, 0, 1$.  So while $B=k[x, 1/x, 1/(x-1)]$ is not free, 
it should in some sense be ``free enough" to be used as a substitute for a free algebra when studying the four-punctured sphere.  Further, since $B$ is a localization of a polynomial ring in one variable, which is free even in the category of non-commutative algebras, it seems reasonable to stay within the commutative category. There,  if an algebra is free (i.e., a polynomial ring), then so is its tensor product with itself (a polynomial ring in twice the number of variables).   So we set  $F=B\otimes B$, which may be identified with $k[x,1/x,1/(x-1),y,1/y, 1/(y-1)]$, with the expectation that this is sufficiently free. There is then an epimorphism $\beta:F \to A= [x, 1/x, 1/(x-1),1/(x-\lambda)]$ sending $x \in F$ to $x \in A$ and $y \in F$ to $\lambda^{-1}x \in A$.  This gives a presentation of $A$ as a quotient of an the algebra $F$  by an ideal $K_{\lambda}$ depending on $\lambda$.  The deformations of the corresponding diagram should capture the dependence of $A$ on $\lambda$. Despite this simplification, to compute the cohomology of the diagram with coefficients in itself requires knowing that of the three algebras $A, F, K_{\lambda}$. The first two are affine but the third is not and we do not know its cohomology.  However, the morphism $\beta:F \to A$ was nothing more than a pair of morphisms 
$k[x, 1/x, 1/(x-1)] \to A$, so with  $B = k[x, 1/x, 1/(x-1)]$. Assuming that $\lambda$ is invertible, let  $f,g:B\to A$ be defined, respectively, by sending $x$ to $x$ and $x$ to $\lambda^{-1}x$,
as pictured in the following diagram which we will denote by $\mathbb{A}$.
\begin{equation}\label{4pt}
 \xymatrix@1{k[x,\frac1x,\frac1{x-1}]=B \hspace{.1in} 
\ar@<+.7ex>[r]^-{ \hspace{.2in}
{f}\hspace{.2in}} \ar@<-.7ex>[r]_-{\hspace{.1in} {g\hspace{.1in}}\hspace{.1in}} & \hspace{.1 in}A= k[x,\frac1x,\frac1{x-1},\frac1{x-\lambda}}].
\end{equation}

\noindent This greatly simplified diagram (\ref{4pt}), whose cohomology is computable, will still detect the dependence of $A$ on $\lambda$. (The next section discusses the problem of viewing a varying presentation as a deformation.)  Geometrically, when $k = \mathcal{C}$, the picture is that of two morphisms of a four-punctured sphere into the three-punctured sphere.

One might hope that the cohomology of the diagram $\mathbb A$ in (\ref{4pt}) would be one dimensional since analytically the  four-punctured sphere has but one parameter. Unfortunately, this is not quite the case. 

\begin{theorem}
The dimension of $H^2(\mathbb{A}, \mathbb{A})$ is infinite. Its generators divide into two classes: The first contains a single (expected) element (corresponding to the number of parameters of the four-punctured sphere). If a natural condition on the regularity of the derivations that are allowed is imposed then the second class, which is infinite, is eliminated.
\end{theorem}
 \noindent\textsc{Proof.}  Both algebras in the diagram (\ref{4pt}) are commutative, so its infinitesimal deformations may be identified with its second cohomology group.  The 2-cocycles here are of the form $\G = (\G^B, \G^A, \G^f, \G^g)$, where $\G^B$ is a 2-cocycle of $B$ with coefficients in itself, similarly for $\G^A$, while $\G^f, \G^g$ are 1-cochains of $B $ with coefficients in $A$, and $f\G^B-\G^Af =\delta \G^f, g\G^B-\G^A g = \delta \G^g$.  Since $B$ and $A$ both have vanishing cohomology in all dimensions greater than one,  $\G$ is cohomologous to a cocycle in which $\G^B = \G^A =0$, which we may now assume, and in this case the cocycle condition is just that $\G^f$ and $\G^g$ must be derivations from $B$ into $A$.  A derivation of any algebra in the diagram is completely determined by its value on $x$ which in turn can be any element of the coefficient module which here is $\mathbb{C}[x, 1/x, 1/(x-1), 1/(x-\lambda)]$.  
An element of this algebra is the sum of its principal parts, i.e, a unique linear combination of $1$ and powers of $x, 1/x, 1/(x-1)$ and $1/(x-\lambda)$.  The 2-cocycle $\G$ can be identified with the pair of elements $(\G^f(x), \G^g(x))$ in $A = k[x, 1/x, 1/(x-1), 1/(x-\lambda)]$.  

To compute the second cohomology of the diagram $\mathbb{A}$ we must compute when a cocycle $\G = (0,0,\G^f, \G^g)$ is the coboundary of a 1-cochain $\g = (\g^B, \g^A, \g^f, \g^g)$. Since $\G^B =\G^A = 0$, the maps $\g^B$ and $\g^A$ must be derivations of $B$ and of $A$, respectively, and $\g^f , \g^g$ (which are just elements of $A$) may be assumed to be zero,  using the asimplicial subcomplex. (Alternatively, only their coboundaries enter but these vanish since the algebras are commutative so their values are irrelevant.)  Thus $H^2(\mathbb{A}, \mathbb{A})$ is the $k$ module of 2-cocycles $(0,0, \G^f, \G^g)$ modulo the submodule of those for which there exist $\g^B \in \operatorname{Der}(B), \g^A \in \operatorname{Der}(A)$ such that simultaneously $\G^f = f\g^B-\g^Af,\, \G^g = g\g^B-\g^Ag$.  The coboundary $\g$ is determined by the elements $\g^B(x) \in B$ and  $\g^A(x) \in A$. Its coboundary corresponds to the pair of elements $f(\g^B(x))-\g^A(fx), g(\g^B(x)) - \g^A(g(x))$. The first is just $\g^B(x) -\g^A(x)$, both summands considered as elements of $A$, while the second is $\g^B(\lambda^{-1}x) - \lambda^{-1}\g^A(x)$.  So $\G$ is a coboundary precisely when there exist elements $b(x) \in k[x, 1/x, 1/(x-1)], a(x) \in k[x, 1/x, 1/(x-1), 1/(x-\lambda)]$ such that (when both are considered as elements of $A= k[x, 1/x, 1/(x-1), 1/(x-\lambda)]$) $b(x) - a(x) = \G^f(x), b(\lambda^{-1}x) - \lambda^{-1}a(x)) = \G^g(x)$. Then $b(x) - \lambda b(\lambda^{-1}x) = \G^f(x) - \lambda\G^g(x)$; if we can solve this then $a(x) = b(x) - \G^f(x)$ is determined. 

The foregoing shows that the cohomology class of the infinitesimal determined by $\G^f, \G^g$ depends only on $\G^f-\lambda\G^g$, so each class contains a representative in which, say, $\G^g = 0$. This could have been anticipated: Intuitively, a derivation of $B$ into $A$ may be viewed as the infinitesimal of a family $g(\hbar)$ of monomorphisms with $g(0) = g$, and similarly with $f$.  Since here $B$ is actually contained in $A$ and every derivation of $B$ into $A$ can be extended to  one of all of $A$, this hypothetical family can be extended to a family (again hypothetical) of automorphisms of $A$. Since only the relation between $f$ and $g$ is important, using the inverses of these automorphisms we could keep $g$ fixed and let only $f$ vary.

The cocycle $\G = (0,0, \G^f, \G^g)$ is a coboundary precisely when we can find $b(x) \in B = k[x, 1/x, 1/(x-1)]$ such that $b(x) -\lambda b(\lambda^{-1}x)= \G^f(x) -\lambda\G^g(x) \in A$. Now $\G^f(x)- \lambda\G^g(x)$ can be an arbitrary element of $A$ so $H^2(\mathbb{A}, \mathbb{A})$ is isomorphic to the $k$-module obtained by taking the quotient of  $A$ as a $k$-module by the submodule consisting of those elements $c(x) \in A$ of the form $b(x) - \lambda b(\lambda^{-1}x)$ (where $\lambda \ne 1$). This is an infinitely generated free module consisting of two natural subfamilies, one with but a single generator. First, we can solve $b((x) -\lambda b(\lambda^{-1}x) = x^N$ for all positive and negative integers $N$ other than $N=1$; this determines a subspace of $H^2(\mathbb{A}, \mathbb{A})$ with but one non-trivial class (up to multiples).  Second,  if $b(x) = 1/(x-1)^N$ then $b(x) -\lambda  b(\lambda^{-1}x) = 1/(x-1)^N - \lambda/ (\lambda^{-1}x-1)^N = 1/(x-1)^N - \lambda^{N+1}/(x-\lambda)^N$ and we can solve $b(x) -\lambda b(\lambda^{-1}x)=c(x)$ if and only if $c(x)$ contains no term in $x$ and the coefficients of $1/(x-1)^N$ and $1/(x-\lambda)^N$ in $c(x)$ are in the ratio $1:-\lambda^{N+1}$.  Therefore, $H^2(\mathbb{A}, \mathbb{A})$ has a basis consisting of classes corresponding to $x$ and all $1/(x-1)^N, N > 0$.

The class corresponding to $x$ is what we have been looking for, but what is the meaning of the   infinite family? The analytic point of view would immediately eliminate it since it arises from allowing the derivations to be singular at $x=1$ and $x= \lambda$. In the geometric picture the tangent vector fields are holomorphic. To follow that picture one should require regularity at $1$ and at $\lambda$; then all the $1/(x-1)^N$ are eliminated. $\Box$\medskip

While the regularity condition eliminates the infinite family, a problem remains:  To show that $A$ actually depends on the parameter $\lambda$ one must show that the diagram (\ref{4pt}) has not only non-trivial infinitesimal deformations but that some at least of these can be integrated into actual deformations. One might suspect that the infinite family consisting of those infinitesimals with singularities are obstructed while the one remaining infinitesimal is integrable, corresponding to the classical variation of $\lambda$. However, it is easy to check that $H^3(\mathbb{A}, \mathbb{A})= 0$ so there are no formal obstructions to any infinitesimal deformations, although over the complex numbers it may be that the infinitesimals in the infinite family have no convergent integrals to actual complex deformations.

To understand the one remaining infinitesimal, simplify matters by assuming that $k = \mathbb{C}$.  One natural integral of the infinitesimal deformation 

$$ \xymatrix@1{k[x,\frac1x,\frac1{(x-1)}]=B \hspace{.1in} 
\ar@<+.7ex>[r]^-{ \hspace{.2in}
{\G^f}\hspace{.2in}} \ar@<-.7ex>[r]_-{\hspace{.1in} {\G^g\hspace{.1in}}} & \hspace{.1 in}A= k[x,\frac1x,\frac1{x-1},\frac1{x-\lambda}}]$$

is then

$$ \xymatrix@C=4pc@1{k[x,\frac1x,\frac1{(x-1)}]=B \hspace{.1in} 
\ar@<+.7ex>[r]^-{ \hspace{.2in}
{\exp(\hbar\G^f)}\hspace{.2in}} \ar@<-.7ex>[r]_-{\hspace{.1in} {\exp(\hbar\G^g)\hspace{.1in}}} & \hspace{.1 in}A= k[x,\frac1x,\frac1{(x-1)},\frac1{(x-\lambda)}}].$$

However, the exponentials must be understood in the sense that the derivations are extended to all of $A$, where the exponentials are now families of automorphisms not of $A$ but of $A[[\hbar]]$ considered as an algebra over $\mathbb{C}[[\hbar]]$, and that using these, $f$ and $g$ have been extended to families of morphisms $f(\hbar), g(\hbar)$ of $B[[\hbar]]$ into $A[[\hbar]]$. The actual picture is

$$ \xymatrix@C=5pc{B[[\hbar]] \hspace{.1in} 
\ar@<+.7ex>[r]^-{ \hspace{.2in}
{\exp(\hbar\G^f)}\hspace{.2in}} \ar@<-.7ex>[r]_-{\hspace{.2in} {\exp(\hbar\G^g)\hspace{.2in}}}& \hspace{.1in}A[[\hbar]]}.$$

Now take $\G^g = 0$ and let $\G^f$ be defined by $\G^f(x) = x$.  In this case we should like to be able to specialize $\hbar$ to an arbitrary element of $\mathbb{C}$ and obtain a deformed diagram of algebras defined over the original ground field.  This deformation would be non-trivial because its infinitesimal was a non-zero class. Specializing $\hbar$ to $\mu \in \mathbb{C}$ does, however, cause a problem since  denoting $\exp(\hbar\G^f)$ by $f_{\mu}$ one has $f_{\mu}x = e^{\mu}x$ and therefore $f_{\mu}(1/(x-1))= 1/(e^{\mu}x-1)$, which is not in the original algebra $A$.  But we have seen that $A[[\hbar]]$ is independent of $\lambda$ and therefore would not change, up to isomorphism, if $\lambda$ were replaced by $\lambda\mu$, so $\mathbb{C}[x, 1/x, 1/(x-1), 1/(x-\lambda)][[\hbar]] \cong \mathbb{C}[x, 1/x, 1/(\mu x-1), 1/(x-\lambda)][[\hbar]]$. The best way to make sense of a specialization of $\hbar$ to $\mu$, therefore, is to view $A[[\hbar]]$ as specializing to the smallest subalgebra over $\mathbb{C}$ containing the images of $B$ under both $g$ and $f_{\mu}$, namely $\mathbb{C}[x, 1/x, 1/(\mu x-1), 1/(x-\lambda)]$.  In some (almost quantum-theoretic) sense any single $\mathbb{C}[x, 1/x, 1/(x-1), 1/(x-\lambda)][[\hbar]]$ thus already contains all $\mathbb{C}[x, 1/x, 1/(x-1), 1/(x-\lambda')]$.

The deformed diagram of $\mathbb{C}[[\hbar]]$ algebras has thus `descended' to an ordinary diagram of algebras over $\mathbb{C}$ and one sees that the class of the cocycle defined by $b(x) = x$ above, which is non-trivial, corresponds to the varying of $\lambda$. So the cohomology of diagrams of algebras has captured what analytically was fairly obvious, namely that $A$ does depend essentially on $\lambda$. What remains unexplained is the meaning of the other, infinite, family of non-trivial classes whose presence seems to be invisible in the analytic picture.  On the other hand, it may be that the diagram (\ref{4pt}), chosen because its cohomology was relatively easy to compute, is not the most appropriate for the problem

A question about the Hochschild-Kostant-Rosenberg theorem arising from the foregoing is whether it can be refined to a similar assertion when certain regularity conditions are imposed on Hochschild cochains.  A more immediate question, however is whether variation of the parameter $\lambda$ is actually a deformation in the classical sense: Is there a fixed basis for $\mathbb{C}[x,1/x, 1/(x-1), 1/(x-\lambda)]$ such that varying $\lambda$ to, say $\lambda(1 +\hbar)$ can be realized as a ``star multiplication'' on that basis?  This is the case, as shown in the next section, but raises the more general question of when variation of a parameter is actually a deformation.

\section{Quotients as deformations and Gr\"obner bases}  To view  a parameterized algebra $A$ defined as a quotient of a fixed ring $R$ by a varying ideal as a classical deformation, one must  display the quotient as a  varying multiplication on a fixed underlying vector space. We seek a criterion for that to be possible.

In the simplest case, suppose that $A$ is an affine ring, i.e., a
finitely generated commutative ring over a field $k$. It is then a
quotient $A = R/I$ with $R$ a polynomial ring $k[x_1,\dots,x_n]$ in
some finite number of variables by an ideal $I$ which itself can be
finitely generated. After choosing a term order on the monomials in
$R$, Buchberger's algorithm produces a unique Gr\"obner basis. The set of (classes of) standard
monomials, i.e., those not in the initial ideal of the Gr\"obner
basis then forms a basis for the quotient algebra $A$. The algorithm
involves only rational operations, so if we have an ideal
$I_{\hbar}$ generated by polynomials which depend rationally on some
formal parameter $\hbar$, then the polynomials in the  basis will
also depend rationally on $\hbar$. Call the poles of these rational functions,  which are finite in number,
the \emph{exceptional values} for the ideal $I_{\hbar}$.  
Except at these values of $\hbar$ in $k$ one will have
the same set of standard monomials and therefore a fixed basis for
 $R/I_{\hbar}$.   With the foregoing notation we therefore have the following.

\begin{theorem}\label{Buchberger}
If zero is not an exceptional value for the ideal $I_{\hbar}$
then  $A_{\hbar}$ may be viewed as a deformation of $A$. $\Box$
\end{theorem}

It would be exceedingly useful if one could extract  from Buchberger's theory some cohomological information which would tell if the deformation is trivial over the original ground field rather than over $k[[\hbar]]$. 

Buchberger's theory is applicable also to certain classes of
non-commutative algebras, but when the reduced Gr\"obner basis is
infinite, even though we may have a fixed complement
to the ideal $I_{\hbar}$ there can be infinitely many exceptional
values in any neighborhood of $\hbar = 0$.  A representation of
$A_{\hbar}$ as a classical deformation may then be impossible. We will see this with the Weyl algebra.

Theorem \ref{Buchberger} can be illustrated with the algebra of the four-punctured sphere.  The algebra $\mathbb{C}[x, 1/x, 1/(x-1), 1/(x-\lambda)]$ is a quotient of the polynomial ring $R =\mathbb{C}[x,y,z,w]$ by the ideal $I_{\lambda}$ generated by $xy-1,
(x-1)z-1, \linebreak[0] (x-\lambda)w-1$. Varying the ideal by varying $\lambda$ can in this case  indeed be viewed as a deformation. Taking as term order total degree with
reverse lexicographic order to break ties (frequently the fastest,
computationally), the resulting Gr\"obner basis (using Maple) is
\{$\lambda zw+z-wz-w, \lambda yw +y -w, wx-\lambda w -1, zy+y-z,
zx-z-1,xy-1$\}. For simplicity we have not divided by $\lambda$ when
that is the leading coefficient. Here $\lambda$ has been treated as
a variable but it is clear from the form of the result that the only
exceptional values of $\lambda$ are $0$ and $1$ (and technically,
$\infty$). Since we start with a fixed value of $\lambda$, the ideal $I_{\hbar}$ above should be viewed as one defined by $\lambda(1+\hbar)$, for which zero is therefore not an exceptional value away from $\lambda = 0,1, \infty$. For every other value of $\lambda$ there is a
neighborhood in which $R/I_{\lambda}$ is indeed a deformation, in the classical
sense, of the algebra of the four-punctured sphere. The relations
$xy-1$ and $(x-1)z-1$ have not changed; the subalgebra generated by
$x,y$ and $z$ is that of the three-punctured sphere, independent of
$\lambda$.  However, the Gr\"obner basis contains a new relation amongst $x,y$
and $z$, namely $zy+y-z$. This is an immediate consequence of the
original defining relations but is necessary: The Gr\"obner basis
shows that the initial ideal is generated by $zw, yw, wx, zy,zx, xy$
so the standard monomials $S$, which span a complementary vector
space $\mathbb{C}S$ to $I(\lambda)$ for all but the exceptional
values of $\lambda$, are those not divisible by any of these. But
this simply says that $S$ consists of all pure powers of either
$x,y,z$ or $w$ and $1$. This has recaptured the decomposition of an element of $k[x, 1/x, 1/(x-1), 1/(x-\lambda)]$ into its principal parts. What remains missing is some indication from the procedure that the quotient algebra depends essentially on the parameter $\lambda$. 

\section{Second example:  The Weyl algebra } Throughout this section the coefficient ring $k$ is assumed to be a field of characteristic zero . The (first) Weyl algebra $W_1 =
k\{x,y\}/(xy-yx-1)$ can be exhibited in various ways as a
deformation of the polynomial ring $k[x,y]$, of which the
simplest is the {\em normal form}: Letting ``$*$'' denote the deformed
multiplication one sets
 $$a*b = ab
+\hbar\partial_xa\partial_yb +\frac{1}{ 2!}\hbar^2
\partial_x^2a\partial_y^2b + \dots =
\mu[\exp(\hbar\partial_x\otimes\partial_y)(a\otimes b)],$$ where on
the right $\mu$ denotes the original multiplication. If $\phi$ and $\psi$ are commuting derivations of an algebra $A$ of characteristic zero then setting $a*b =
\mu\exp(\hbar\phi\otimes\psi)(a\otimes b)$ defines a deformation of
$A$; more generally, given commuting derivations
$\phi_1,\dots, \phi_r, \psi_1,\dots,\psi_r$ one can replace
$\phi\otimes\psi$ by $\sum \phi_i \otimes \psi_i$,  \cite{G:DefI, G:DefIV}. In particular,
$\phi\otimes\psi$ can be replaced by $\phi \wedge \psi = (\phi
\otimes \psi - \psi \otimes \phi)/2$ to give a ``Moyal deformation''
\cite{Moyal} (but the idea is already implicit in Groenewold
\cite{Groenewold}). This gives $[x,y]_* = \hbar$; taking $\hbar=1$
gives the Weyl algebra.  In most physical applications the ground ring is $\mathbb{R }$ or $\mathbb{C}$, the algebra $A$ being deformed is one of functions, and  the deformed product is of the form
$a\ast b = ab +\hbar F_1(a,b) + \hbar^2F_2(a,b) + \dots$, where the $F_i$ are (bi)differential operators of increasing orders. One hopes for convergence at least for sufficiently small real or complex values of $\hbar$.  Here we are dealing with polynomials so there is no problem of convergence since the series for any particular product will terminate.  (A polynomial will ultimately be annihilated by a sufficiently high-order differential operator.)   Note that if both $k[x,y]$ and $W_1$ are viewed as singly graded with
 $\deg x =1$ and $\deg y = -1$ then the Moyal deformation has preserved this grading.  (In deformation quantization, the deformation parameter is actually $1/\hbar$.)

The cohomology of $k[x,y]$ is non-trivial, a special case of the theorem of Hochschild-Kostant-Rosenberg \cite{HKR} noted in the Introduction. However that of $W_1$ vanishes in all positive dimensions, a classic theorem of Sridharan \cite{Sridharan:Filtered}, so $W_1$ is absolutely rigid. Nevertheless, one can
introduce a parameter by defining $W_q = k\{x,y\}/(qxy-yx-1)$, and $W_q$ is not isomorphic to $W_1$ for $q \ne 1$. Classical deformation theory can not recognize the passage from $W_1$ to $W_q$ but diagram cohomology can. The  cohomology ring  $H^*(W_q,W_q)$ is not trivial; it is computed in  \cite{GerstGiaq:WeylCohom} which also gives a very short proof of Sridharan's theorem, both using the tools of deformation theory.  In particular, there is a non-trivial $2$-cocycle and classical deformation theory does recognize the variation of $W_q$ with $q$ once $q$ is different from one, but all classes of positive dimension are fragile at $q=1$ since those of $W_1$ all vanish. 

The algebra $W_q$, like $W_1$, is well-known in quantum mechanics, usually in the form of an algebra of operators generated by multiplication by $x$ and $d/\!dx$, since as operators one has $x(d/\!dx) -(d/\!dx) x = -1$, so $-d/\!dx$ can take the place of $y$. We present first some elementary facts which hold for all $q$ including $q=1$ unless stated otherwise.  First, as mentioned, there is a $\mathbb{Z}$ grading which survives from the bigrading of $k\{x,y\}$: define the degree of a monomial 
$x^{i_1}y^{j_1}x^{i_2}y^{j_2}\cdots x^{i_r}y^{j_r}$ to be  $\sum i_k - \sum j_k$ and define an element to be homogeneous of degree $m$ if it is a linear combination of monomials of degree $m$. This is well-defined and independent of the way the element is 
written, because this degree is well-defined in $k\{x,y\}$  and the equation defining $W_q$ is homogeneous of degree zero. Henceforth, by the degree of an element of $W_q$ we will mean only this degree.
Every element of $W_q$ can be written uniquely as a {\em pseudopolynomial}, i.e., a linear combination (with coefficients in $k$) of monomials of the form $x^iy^j$.

 This is a 
special case of the following, whose proof is but omitted.
\begin{theorem}
 Suppose that an algebra $A$ has an exhaustive ascending filtration, $F_0A \subset F_1 A \subset F_2A \subset \dots; \cup F_iA = A$ (where $F_iA\cdot F_jA 
\subset F_{i+j}A$) with the property that the associated graded ring $\operatorname{gr}A$ has a set of generators $\bar S$ such that for every $\bar a, \bar b \in \bar S$ one has $\bar b\bar a = \lambda \bar a\bar b$ for some $\lambda \in k$. Then $A$ has a set of generators $S$ (representatives of the 
elements of $\bar S$) such that ordering arbitrarily the elements of $S = \{x_1,\dots x_n\}$ (which for ease of notation is assumed finite), every element of $A$ can be rewritten as a linear combination of 
monomials $x_1^{i_1}x_2^{i_2}\cdots x_n^{i_n}$. Moreover, if the representations in $\operatorname{gr}A$ are unique then so are the latter. $\Box$
\end{theorem} \medskip

In $W_q$, take $F_iW_q$ to be the span of all elements $a$ with $|\deg a| \le i$. 

\begin{lemma}
The elements of $W_q$ of degree zero form a commutative subalgebra. 
\end{lemma}

\noindent\textsc{Proof.}  That they form a subalgebra is evident. A simple induction shows that writing $[n]_q = 1+q+q^2+\cdots q^{n-1}$ the defining equation $qxy - yx = 1$ implies 
\begin{equation}
q^nxy^n-y^nx = n_qy^{n-1}.
\label{BasicIdentity}
\end{equation}
With this, another simple induction also shows that
\begin{equation*}
x^ny^n = (xy+[n-1]_q)(xy+[n-2]_q)\cdots(xy+[1]_q)xy.
\end{equation*}
The right side is a polynomial in the single element $xy$ where when expanded the coefficients of the powers of $xy$ are (up to sign, depending on choice of definition) the $q$-Stirling numbers of the first kind; these are amongst those introduced by Carlitz \cite[p. 129]{Carlitz:Abelian}.  Conversely, one can rewrite the powers of $xy$ as linear combinations of the $x^ny^n$ in the usual way using the $q$-Stirling numbers of the second kind; it follows that the subalgebra of $W_q$ spanned by elements of degree zero is just the polynomial subalgebra generated by the single element $xy$.  $\Box$\medskip

The $q$-Weyl algebra and the $q$-Stirling numbers (both of the first and second kind) have applications in quantum theory in the study of the $q$-harmonic oscillator and $q$-creation and annihilation operators, cf. eg. \cite{MendezRodriguez:q-Stirling}, \cite{Parasarathy:q-Fermionic}. 

Since $[n]_q$ vanishes when $q$ is a primitive $n$th root of unity, equation (\ref{BasicIdentity}) shows that when $q$ is such a root then $x^n$, and similarly $y^n$, are central in $W_q$. It is easy to check that these generate the center, which is then just the polynomial ring $k[x^n,y^n]$.  By contrast, the center of $W_1$ is reduced to the identity. This shows that any ordinary moduli space for the $q$-Weyl algebras must be very badly behaved near $q=1$ and every root of unity, something which will manifest itself shortly in another way.

As shown in e.g. \cite{GerstGiaq:WeylCohom},\cite{Witherspoon:Hochschild}, $W_q$ has vanishing cohomology in all dimensions greater than two. When $q$ is not a root of unity, which we will henceforth assume, $H^1(W_q, W_q)$ is generated in dimension one by the class of a single cocycle denoted $x\partial_x-y\partial_y$ and $H^2(W_q, W_q)$ is generated also by the class of a single cocycle. To understand these we must first consider the cohomology of the algebra 
 $W_{\text{qp}}$  obtained by imposing on the non-commutative polynomial ring $k\{x,y\}$  (the free unital associative algebra) in two variables the relation $qxy=yx$; this is the algebra of the quantum plane.  This algebra is a deformation of $k[x,y]$ induced by the infinitesimal deformation  $x\partial_x \wedge y\partial_y$ of the latter.  When $q$ is not a root of unity, $H^2(W_{\text{qp}}, W_{\text{qp}})$ is spanned by the classes of exactly two cocycles.  The first is the lifting of the class of the cocycle $\partial_x \wedge \partial_y$. Note that this is the class which induces the deformation of $k[x,y]$ to the Weyl algebra; while it lifts to the Weyl algebra it becomes a coboundary there, as do all two-dimensional classes.  (We have availed ourselves here of the abuse of language mentioned near the beginning, since in the deformed algebra these classes are only $\hbar$ torsion classes until $\hbar^{-1}$ is adjoined.) However $\partial_x \wedge \partial_y$ lifts non-trivially in  $W_{\text{qp}}$ to a cocycle whose explicit form is not needed but which will be denoted by $z_{\text{qp}}$. The $q$-Weyl algebra with deformation parameter a variable $\hbar$, given by  $k\{x,y\}/(xy-qyx-\hbar)$ and denoted for the moment by $W_q(\hbar)$ (so $W_q(1) = W_q$), is a jump deformation of $W_{\text{qp}}$;  it is isomorphic, except for extension of coefficients, to $W_q$.  (This justifies the abuse, since when $\hbar$ is specialized to $1$ it is certainly invertible.) The cocycle $z_{\text{qp}}$ is actually the infinitesimal inducing the jump deformation of $W_{\text{qp}}$ to $W_q$; it lifts to a coboundary (as does the infinitesimal of any jump deformation).  The second class is that of $x\partial_x\wedge y\partial_y$. This is a cocycle of $k[x,y]$ which
 in fact continues to be well-defined in $W_{\text{qp}}$ (since the same is true of $x\partial_x$ and $y\partial_y$);  it lifts to itself.

Setting $q=1$ in $W_q$ gives the Weyl algebra, whose cohomology vanishes in positive dimensions, so one might expect that the same would be the case for $W_q$ (which initially appears to be a deformation of $W_1$), but this is not the case. The $2$-cocycle $x\partial_x\wedge y\partial_y$ of $W_{\text{qp}}$ lifts non-trivially to a $2$-cocycle $\hat z$ of $W_q$. Its explicit form again is not important; it is the infinitesimal of the deformation of $W_q$ to $W_{q'}$ with $q'$ near $q$. It follows that for $q \ne 1$, $W_q$ is not absolutely rigid and that its deformation to $W_{q'}$ can be detected by the classic deformation theory. (The second cohomology of $W_q$ has more non-trivial classes when $q$ is a root of unity; for details, cf. \cite{GerstGiaq:WeylCohom}.  For simplicity we will continue to assume that $q$ is not a root of unity.)

The algebras $W_1$ and $W_q$ share the property that every element can be written uniquely as a pseudopolynomial, so their underlying vector spaces can be identified. It would seem, then, that one could exhibit $W_q$ with generic $q$ as a deformation of $W_1$, but as $W_1$ is absolutely rigid, that is in fact impossible. Nevertheless, we have the following.

\begin{theorem}
Setting $q=1+\hbar$,
 $W_1[[\hbar]]$ and $W_q[[\hbar]]$ are isomorphic as algebras over $k[[\hbar]]$.  All cohomology classes of  $W_q[[\hbar]]$, like those of $W_1$, vanish in positive dimensions; all of its higher dimensional classes are therefore fragile and the resilient cohomology of $W_q$ is reduced just to the constants.  
\end{theorem}
\noindent\textsc{Proof.} The second assertion follows from the first. For the first we show that
in fact, with $q=1+\hbar$, $W_1[[\hbar]]$ and $W_q[[\hbar]]$ are not merely isomorphic but identical. (We could, as later, set $q= e^{\hbar}$; in either case, $q(0) = 1$.) To prove this it is sufficient to exhibit an element $z \in W_q[[\hbar]]$ with $z \equiv y \pmod{\hbar}$, i.e., of the form $z = y+\hbar\eta_1 +\hbar^2\eta_2 + \dots,  \eta_i \in
W_q$ such that $xz-zx = 1$.

One can solve for the $\eta_i$ recursively. Reducing the relations
$$[x, y+\hbar\eta_1 + \hbar^2\eta_2 + \dots], \quad [x,y] = 1-\hbar xy$$ modulo
$\hbar^2$ gives $[x,y] + \hbar[x,\eta_1] \equiv 1 \pmod \hbar^2 $ or
$1-\hbar xy + \hbar[x,\eta_1] \equiv 1 \pmod \hbar^2$, whence
$[x,\eta_1] \equiv xy \pmod \hbar$. Here there is some choice. Taking
 $\eta_1 = (1/2) xy^2$ (the simplest choice), reduction modulo
$\hbar^3$ gives $[x,y + (\hbar/2) xy^2 + \hbar^2\eta_2] \equiv 1
\pmod {\hbar^3}$ or $[x,y] + (\hbar/2) x\{[x,y]y + y[x,y]\} + \hbar^2
[x,\eta_2] \equiv 1 \pmod{\hbar^3}$. Using again that $[x,y] = 1-\hbar xy$ we have $-\hbar xy + (\hbar/2)x((1-\hbar xy)y + y(1-\hbar
xy)) + \hbar[x,\eta_2] \equiv 0 \pmod {\hbar^2}$. Therefore $[x, \eta_2] =
\frac12(x^2y^2 + xyxy) \bmod \hbar$. Since this is a congruence
modulo $\hbar$, the computations may now take place in $W_1$.
 In the previous congruence we should therefore replace
the occurrence of $yx$ by $xy-1$. Again making the simplest
choice, we can take $\eta_2 = \frac13 x^2y^3 - \frac14 xy^2$.
Continuing in this way one finds the following recursion for getting
$\eta_{r+1}$ from $\eta_r$. First, apply $xy\partial_y$, then
reduce the result to a pseudopolynomial in $W_1$ using $[x,y] =
1$, and finally integrate with respect to $y$. The result is this:
If $p(x,y)$ is any pseudopolynomial in $x$ and $y$, define
$\hat{p}$ to be the result of replacing $y^n$ everywhere by $
(n/(n+1)) xy^{n+1} - ((n-1)/2) y^n$. Then $\eta_{r+1} =
x\hat{\eta}_r$.

With this one can compute compute 
\begin{equation*}
\begin{split}
z &= y+\hbar\eta_1 +
\hbar^2\eta_2 + \dots\\ & = y + \hbar(\frac12 xy^2) +\hbar^2(\frac13
x^2y^3- \frac14 xy^2) + \hbar^3(\frac14x^3y^4 - \frac12x^2y^3 +\frac14 xy^2) +
\dots
\end{split}
\end{equation*}
This exhibits the desired $z$.$\Box$ \medskip

In the preceding proof, the coefficient $\eta_1$ of $\hbar$ in $z$, namely $(1/2)xy^2$ can be viewed as an infinitesimal deformation, as will be seen below.

Only monomials of the form $x^my^{m+1}, m = 0, 1, 2 \dots$ appear in
$z$ above, which suggests reordering the series into one of the form $y+a_1xy^2 +
a_2x^2y^3 + \dots$, the coefficients $a_i$ now being power series in
$\hbar$. As the relation $qxy-yx = 1 $ implies both
\begin{equation}\label{recursion}
q^nx^ny-yx^n =\frac{q^n-1}{q-1}x^{n-1} \quad \text{\ and\ }  \quad
q^nxy^n -y^nx = \frac{q^n-1}{q-1}y^{n-1}
\end{equation}
 one obtains by recursion that 
$$z = y + a_1(\hbar)xy^2 + a_2(\hbar)x^2y^3 + \dots, \quad \mathrm{\
where\ \ } a_r(\hbar) = \frac{\hbar^{r+1}}{(1+\hbar)^{r+1}-1}.$$
In the power series representation of $z$ there could be no value of $\hbar$ for which the
series converges, for each $a_r(\hbar)$ has 
a pole wherever $\hbar$ has the form $\omega\! -\! 1$ where $\omega$
is an $(r+1)$st root of unity and every neighborhood of 0 in
$\mathbb{C}$ contains infinitely many of these. The extraordinary behavior of $W_q$ at the roots of unity is thus reflected in the deformation theory; those $W_q$ with
$q$ a root of unity are in some sense `unreachable' by deformation from $W_1$.  As $q=1$ is a limit point of the roots of unity, which are in fact dense on the unit circle, it is clear that no sort of formally analytic deformation (as in the classical theory) from $W_1$ to $W_q$ is possible.
Nevertheless, $z$ can actually be evaluated for any complex
$\hbar$ with $1+\hbar$ not a root of unity.

There is another isomorphism $W_1[[\hbar]] \cong
W_{1+\hbar}[[\hbar]]$, this time starting with $W_1$, given by the second author and Zhang,
\cite{GZ:QuantWeyl}. Writing $D$ for $d/\!dx$ one has (as operators) the classic relation
$xD - Dx = -1$, so $x$ and $-D$ satisfy the same relation as do
$x$ and $y$ in $W_1$. The element
$$\zeta =\frac{(1/x)(e^{\hbar xD}-1)}{e^{\hbar} -1}$$ then has the property
that $\zeta x  = e^{\hbar}x\zeta + 1$, for note that $$\zeta x =
x^{-1}\frac{e^{\hbar xD}-1}{e^{\hbar}-1}x = (e^{\hbar
x(x^{-1}Dx)}-1)/(e^{\hbar}-1)\quad \mathrm{and} \quad x^{-1}Dx =
D+ x^{-1},$$ from which the equation follows. Therefore, setting
$y_{\hbar} =-x^{-1}(e^{-\hbar xy} -1)/(e^{\hbar} - 1)$ one has
$e^{\hbar}xy_{\hbar} - y_{\hbar} x = 1$. Here $e^{\hbar}$ takes the place of $q$ and the expression becomes infinite, as before, when $q$ is a root of unity. 

To capture homologically the passage from $W_1$ to $W_q$ consider the diagram 
$$ \begin{CD} B=k[x]@>f>>W_1 @<g<< k[y] = C\end{CD} $$ 
where $f$ and $g$ are the inclusion maps. We will denote it by $\mathbb{W}_1$.

\begin{theorem}
A basis for the infinitesimal deformations of the diagram $\mathbb{W}_1$ can be identified with the set of monomials $x^iy^j \in W_1$ with either $i >0, j=0$ or $i>0, j\ge 2$ (i.e., omit 1, powers of $y$ alone, and monomials of the form $x^ry$, all $r$).
\end{theorem}
\noindent\textsc{Proof.}
Since all the algebras in the diagram  $\mathbb{W}_1$ are absolutely rigid, every two cocycle of that diagram is cohomologous to one of the form $\G = (0,0,0,\G^f, \G^g)$, the only restriction on which is that $\G^f, \G^g$ be derivations of $k[x]$ into $W_1$ and of $k[y]$ into $W_1$, respectively. Such a cocycle is a coboundary precisely when there are derivations 
$\g^B, \g^C, \g^{W_1}$ of $B,C,W_1$, respectively, such that $\G^f = f\g^B-\g^{W_1}f, \,\G^g = g\g^C - \g^{W_1}g$. Every derivation of $W_1$ is inner, so there is $\chi \in W_1$ such that $\G^{W_1} = \operatorname{ad}{\chi}$. Then $\G^f(x) = \alpha(x) - [\chi, x]$ for some polynomial $\alpha$ and similarly $\G^g(y) = \beta(y) -[\chi, y]$ for some polynomial $\beta$. The first equality can always be satisfied by choosing $\xi$ suitably, so every two-cocycle $\G$ is cohomologous to one of the form $(0,0,0,0, \G^g)$ where now we can only replace $\G^g$ by $\G^g -(\g^Cg-\operatorname{ad}{\xi})$ with $\xi$ having the property that $[\xi, x]$ is a polynomial in $x$. Therefore, up to a constant, $\xi$ must be a linear combination of monomials of the form $x^ry$. Finally, $\G^g(y)$ is an arbitrary polynomial in $y$, so finally the infinitesimal deformations of $\mathbb{W}_1$ can be identified with the linear combinations of those monomials of the form $x^iy^j$ omitting those which are just powers of $y$ (including $1$) or of the form $x^ry$ for some $r$. $\Box$\medskip

The interpretation of the pseudopolynomials which appear in the theorem as infinitesimal deformations is clear: if $\mu$ is one of them then we should like to replace $y$ by some $y'$ which is a power series in $\hbar$ beginning $y' = y+\hbar\mu +\hbar^2\mu_2 +\hbar^3\mu_3 + \dots$. This is exactly what was done at the beginning of this example where we had 
$\mu = \frac12 xy^2,\quad \mu_2 =\frac13x^2y^3,\quad  \mu_3=\frac14x^3y^4 - \frac12x^2y^3 +\frac14 xy^2, \dots.$  
It is easy to check that $H^3(\mathbb{W}_1, \mathbb{W}_1)= 0$ so there are no obstructions to any infinitesimal. This is also evident from the fact that we could choose any power series in $\hbar$ for $y'$ and obtain a deformation, and as long as the leading term did not represent a coboundary it would be non-trivial.  The series exhibited was the one that produced $W_q$ where $q$ was $1+\hbar$ or $e^{\hbar}$.  The defining equation for the resulting algebra in this case can be written in a simple closed form but it is not clear if this is always be possible.  Recall that $W_1$ inherited a single grading from $k[x,y]$ in which the degree of $x^iy^j$ was $i-j$. To preserve it would require that $\mu$ and all the $\mu_i$ above have degree $-1$. The simplest such $\mu$ is $xy^2$.  The diagram cohomology indicates that the simplest way to introduce a parameter into the Weyl algebra while preserving this grading is to pass to the $q$-Weyl algebra.
The question of the existence of a simple closed form for defining equations is related to that of descent,  as in the case of the four-punctured sphere,  but is more difficult here.  

This concludes our examples. The diagrams attached here to the algebras of the examples have done what was intended, namely, \emph{they have exhibited a dependence on parameters which could not be detected by the classical deformation theory of a single algebra.} The deformations of these diagrams are the ``deformations associated with rigid algebras'' of the title.  

While some justification has been given for the construction of these diagrams, it is still largely {\itshape ad hoc} and a unifying theory is needed. However, 
both examples demonstrate  that to understand the deformation of algebras one must consider not just single algebras but diagrams. Also, the example of the Weyl algebra shows that classical deformation is not an inverse operation to specialization (and can not be whenever there are fragile classes), and the moduli space of something seemingly tame, like the Weyl algebra, can be exceedingly complicated. 

\nocite{HazGer}

\end{document}